\documentclass[preprint,authoryear,12pt]{elsarticle}
\usepackage{amssymb}
\usepackage{color}
\begin{document}
\renewcommand*{\arraystretch}{1.3}

\begin{center}
{\Large\bf Constrained discounted  Markov decision processes
with Borel state spaces}
\end{center}
\begin{center}
{\bf \large  Eugene A. Feinberg$^{a}$, Anna Ja\'skiewicz$^{b}$, Andrzej S. Nowak$^{c}$ }
\end{center}
\noindent$^{a}$ Department of Applied Mathematics and Statistics,   Stony Brook University,
             Stony Brook, NY, 11794-3600, USA,
{\footnotesize  email: {\it eugene.feinberg@stonybrook.edu}}\\
\noindent$^{b}$Faculty of Pure and Applied Mathematics, Wroc{\l}aw University of Science and Technology,
             Wybrze\.ze Wyspia\'nskiego 27, 50-370 Wroc\l aw, Poland,
{\footnotesize  email: {\it anna.jaskiewicz@pwr.edu.pl}}\\
\noindent$^{c}$Faculty  of
Mathematics, Computer Science
and Econometrics,  University of Zielona G\'ora,
Szafrana 4$a$, 65-516 Zielona G\'ora, Poland,
{\footnotesize  email: {\it a.nowak@wmie.uz.zgora.pl }}\\
\begin{center}
 \today
\end{center}


\noindent
{\bf Abstract. }
We study  discrete-time discounted  constrained Markov decision processes (CMDPs) with
Borel state and action spaces.
These CMDPs satisfy either weak (W) continuity conditions, that is, the transition probability is weakly
continuous  and the reward function is upper semicontinuous in state-action pairs, or setwise (S) continuity
conditions, that is, the transition probability is setwise continuous and the reward function is upper semicontinuous in actions.
Our  main goal is to study models with unbounded reward functions, which
are often encountered in applications, e.g., in consumption/investment problems.
We provide some general assumptions under which
 the optimization problems in  CMDPs are solvable in the class of randomized   stationary policies and in the class of
 chattering policies introduced in this paper.
If the initial distribution and transition probabilities are atomless, then using a general
``purification result'' of Feinberg and Piunovskiy we show the existence  of a deterministic (stationary) optimal policy.
Our main results are illustrated by examples.

\section{Introduction}

In this paper, we study discrete-time constrained discounted
Markov decision processes on
Borel spaces with unbounded reward functions.
A common feature of almost all previous studies of CMDPs with Borel state spaces is the assumption that the transition
probabilities  (reward functions)  are jointly weakly continuous (upper  semicontinuous) in state-action pairs.
This condition played a significant role  in the convex analytic approach in exploiting properties
of the sets of occupancy measures, see \cite{gmmor,mp,pbook}, \cite{dufour}  and \cite{z}.
The same continuity assumptions were also applied  by \cite{fp1,fp2},
who examined Markov deterministic optimal policies  in non-homogeneous Markov decision processes (MDPs). 
\cite{fr} studied constrained absorbing (in particular, discounted) 
MDPs with both weakly and setwise continuous  transition probabilities.  
MDPs with setwise continuous transition probabilities and multiple criteria were also considered in \cite{fp18}, because their results 
for such a model with finite action sets were used there to prove that, 
if an MDP is atomless (that is, initial and transition probabilities are atomless),  
for an arbitrary policy there exists a non-randomized stationary policy with the same vector of the expected total costs.  
This result is repoerted in \cite{fp18} for uniformly absorbing (in particular, discounted) atomless MDPs 
without any continuity assumptions on costs and transition probabilities. 

This paper studies discounted  CMDPs with   Borel state and action spaces and unbounded reward functions,
where the transition probabilities are weakly or setwise continuous.
A basic ingredient in our approach is the equivalence of the weak topology and the
so-called $ws^\infty$-topology (introduced by \cite{schal1})
on  the space of probability measures induced by policies.
This fact was first observed by \cite{now} and then generalized by \cite{b},
who dealt with the action sets depending on partial histories of the underlying process.
The aforementioned results  allow us to work  in the space of strategic measures and
to prove the existence of a solution,  say $\pi^*,$  to a discounted CMDP
(with unbounded from above and below reward functions) in the class of history dependent policies.
Furthermore, making use of certain lemmas on occupancy measures given under various assumptions
by \cite{bor}, \cite{pbook} and \cite{fr},
we show that $\pi^*$ can be replaced by a randomized stationary optimal policy $\varphi^*.$
In the next step, we assume that the initial distribution and the transition probability
are atomless. Applying a recent  theorem on the performance sets in vector-valued MDPs
due to \cite{fp18}  we ``purify'' $\varphi^*$  and obtain
a non-randomized stationary optimal policy.

In some sense this paper complements the results of \cite{fr},
who studied absorbing CMDPs with rewards bounded from above.
Absorbing MDPs are more general than discounted MDPs. In certain cases discounted
MDPs with unbounded rewards can be transformed to discounted models with bounded
and bounded from above rewards (see Chapter 5 in \cite{vdw}). Then,
the results of \cite{fr} can be also applied to those MDPs.
In this paper, we use assumptions under which the standard transformations do not work. 
This paper is also relevant to the work of \cite{dufour},
who considered unbounded costs with weakly continuous transitions
and some finer  topology on the space of strategic measures.
By imposing uniform integrability  conditions for the positive parts of the reward functions in our model,
we use the standard weak topology on
the space of strategic measures induced by policies. We show that the assumptions imposed by \cite{dufour}
imply our assumptions and that the topology  introduced there is the standard weak topology 
under assumptions of that paper.  

The related literature on the models of CMDPs with finite or countable state spaces is  quite  
large.
The reader is referred to \cite{altmmor,alt, borsiam,fs1,fs2,fs3,kal,pbook} and  references cited therein.
The significance of CMDPs for various
applications  is very-well documented.
The models with discrete state spaces are often used to study admission or flow control problems in queueing networks,
see  \cite{alt,hosp,l,pbook,sen} and \cite{vakil} among others.
The models with uncountable Borel state spaces, on the other hand, are fundamental for  inventory systems,
consumption/investment problems, and some issues in mathematical finance or insurance, see \cite{fp1,fp2,mp,pbook,z}.
The methods which are used in the study of CMDPs are based on  linear programming approach
(see, e.g., \cite{alt,kal,gmmor,dufour}), convex analysis
(see, e.g., \cite{borsiam,fr,fs1,fs2,mp,pbook,z}, Lagrange multipliers (see \cite{altmmor,alt,pbook}),
and dynamic programming techniques (see, e.g.,
\cite{cb,cf,pm}).

%%%%%%%%%%%%%%%%%%%%%%%%%%%%%%%%%%%%%%%%%%%%%%%%%%%%%%%%%

This paper is organized as follows. Section 2 describes the model.
Section 3 contains results on the existence of a solution
in the class of randomized stationary policies. In Section 4,
we study some special classes of randomized stationary policies 
called ``chattering'' and present a  result on the existence
of deterministic  stationary optimal  policies in a class of CMDPs with atomless transitions.
Section 5 contains three illustrative examples related to models developed 
in operations research or economics.
Finally, in Sections 6 and 7 we give additional comments on  our main assumptions.

 \section{The model}

Let $\mathbb{N}$ be the set of all positive integers and $\mathbb{R}$ be the set of all real numbers.
Moreover, put $\mathbb{R}_-:=\mathbb{R}\cup\{-\infty\}.$
By a Borel space, say $Y$,  we mean a non-empty Borel subset of a complete separable metric space.
 Let ${\cal B}(Y)$  denote the $\sigma$-algebra of all Borel subsets of $Y$ and $\Pr(Y)$ be the space
of all probability measures on ${\cal B}(Y)$ endowed with the weak topology. This is the coarsest topology on $\Pr(Y)$
in which the functional  $\nu\to \int_Y u d\nu$ is upper semicontinuous for every
bounded above upper semicontinuous function $u:Y\to\mathbb{R}_-.$  For any $B\in {\cal B}(Y),$
by $1_B(\cdot)$ we denote the indicator function of the set $B.$

By a Borel measurable transition probability from $Y$
to a Borel space $Z$ we mean a function  $\gamma: {\cal B}(Z ) \times Y\to  [0, 1]$ such that, for each $B\in{\cal B}(Z),$
$\gamma(B,\cdot)$ is
a Borel measurable function on $Y$ and $\gamma(\cdot,y) \in \Pr(Z)$ for each $y\in Y.$
We shall write $\gamma(B|y)$ for $\gamma(B,y).$

Let $I:=\{1,\ldots,m\}$ and $I_0=\{0\}\cup I$, where $m\in\mathbb{N}.$

A {\it constrained Markov decision process} (CMDP)  is characterized by the objects:
$S,$ $A,$ $\mathbb{K},$ $p,$ $\mu,$ $r_0,r_1,...,r_m,$ $\beta$  with the following meanings.
\begin{enumerate}
\item[(i)] $S$ is a Borel {\it state space}.
\item[(ii)] $A$ is the {\it action space} and is also assumed to be a Borel space.
\item[(iii)] $A(s)$ is a non-empty
{\it set of actions} available in state $s\in S.$
It is assumed that the set
$$\mathbb{K}  := \{(s, a): \ a \in A(s), s\in S\}$$
is Borel in $S\times A.$
\item[(iv)] $p$ is  a {\it transition probability} from $\mathbb{K}$ to $S.$
\item[(v)] $\mu \in \Pr(S)$ is an {\it initial distribution.}
\item[(vi)] $r_0: \mathbb{K}\to \mathbb{R}_-$ is a Borel measurable {\it stage reward function.}
\item[(vii)] $r_i:\mathbb{K}\to \mathbb{R}_-$, $i\in I$, are Borel  measurable {\it sources of constraints}.
\item[(viii)] $\beta\in (0,1)$ is a {\it discount factor}.
\end{enumerate}

Let $H_n$ be the space of all feasible histories up to the $n$-step, i.e.,
$$H_n=\mathbb{K}^{n-1}\times S\quad \mbox{for }n\ge 2\quad \mbox{and}\quad  H_1=S.$$

An element of $H_{n}$ is called a {\it partial history} of the process
and is of the form
$$h_1=s_1\quad\mbox{and}\quad h_{n} :=
(s_{1},a_{1},\cdots,s_{n-1},a_{n-1},s_{n}), \quad n\ge 2.$$
A  {\it policy} for the decision maker is a sequence $\pi=(\pi_n)$ of transition probabilities
$\pi_n$ from $H_n$ into $A$ such that $\pi_n(A(s_n)|h_n)=1$ for
all $h_n\in H_n$ and $n\in\mathbb{N}.$
The {\it set of all policies} is denoted by $\Pi$.
By $\Phi$ we denote the set of all Borel measurable mappings  $\varphi: S\to\Pr(A)$
such that $\varphi(A(s))(s)=1$ for each  $s\in S.$ Every $\varphi \in \Phi$
induces a transition probability $\varphi(\cdot|s)= \varphi(s)(\cdot)$ from $S$ to $A.$
{\it Markov policy} is a sequence $\pi=(\varphi_n),$ where each $\varphi_n\in \Phi$ for $n\in\mathbb{N}.$
The set of all Markov policies is denoted by $\Pi_M$.
A {\it stationary policy} is a constant sequence $\pi=(\varphi,\varphi,\ldots),$ where $\varphi\in \Phi,$ 
and is identified with $\varphi.$ Therefore, the set of  all stationary  policies   will be also denoted by $\Phi.$
If the support of each  measure $\varphi_n(s)(\cdot)$ is a single point  for every $s\in S$,
then $\pi =(\varphi_n)$ is called {\it non-randomized}
or {\it deterministic} Markov (stationary) policy.
If every set $A(s)$ is compact, then by the Arsenin-Kunugui theorem (see Theorem
18.18 in \cite{k}) the correspondence
$s\to A(s)$ admits a Borel measurable selector, that is,
a mapping $f: S\to A$ such that $f(s)\in A(s)$ for every $s\in S.$
We use $F$ to denote both the set of all such selectors and the set of all  deterministic stationary policies. 

Let $((S\times A)^{\infty},\cal T)$ be the
measurable space, where $\cal T$ denotes
the corresponding product $\sigma$-algebra. Due to
the theorem of Ionescu-Tulcea (see  Proposition V.1.1
in \cite{neveu}), for each
policy  $\pi \in \Pi$
there exists a unique probability measure  $P_{\mu}^{\pi}$ on $\cal T$
such that  for all $D \in {\cal B}(A)$,
$B \in {\cal B} (S)$ and $h_{n} =
(s_{1},a_{1},\ldots,s_{n-1},a_{n-1},s_{n})$ in $H_{n}$, $n \in\mathbb{N}$,
$$P_{\mu}^{\pi}(s_{1}\in B) =\mu(B),$$
$$P_{\mu}^{\pi}(a_{n} \in D| h_{n}) = \pi_{n}(D|h_{n}),$$
$$P_{\mu}^{\pi}(s_{n+1} \in B|h_{n},a_{n}) =
p(B|s_{n},a_{n}).$$
By $E_{\mu}^{\pi}$ we denote the expectation operator
with respect to the probability measure  $P_{\mu}^{\pi}.$\\
Define
\begin{equation}
\label{CP}
J_i(\pi):=E_{\mu}^\pi\left(\sum_{n=1}^\infty \beta^{n-1}r_i(s_n,a_n)\right), \quad i\in I_0.
\end{equation}
Below, we formulate assumptions under which the functionals $J_i(\pi),$ $i\in I_0$ are well-defined for every $\beta\in (0,1).$
\\

\noindent
{\bf Problem Statement.} Let $d_1,\ldots, d_m$  be fixed real numbers.
Consider the following control  problem:
$$\mbox{(CP)}\qquad \mbox{maximize}\quad J_0(\pi)$$
 $$ \mbox{subject to}\quad
 J_i(\pi)\ge d_i, \quad i\in I.
$$
In the sequel, we shall tacitly  assume that problem (CP) is non-trivial, i.e., the set of feasible policies is non-empty and
there exists a feasible policy $\pi$ such that
 $J_0(\pi)>-\infty.$

We now make  the following assumptions:\\

\begin{enumerate}
\item[(A1)] There exists a Borel measurable function $w:S\to[1,\infty)$
such that
$r_i(s,a)\le w(s)$ for each $(s,a)\in \mathbb{K}$ and
for all $i\in I_0$.
\end{enumerate}

Furthermore, we shall also consider the following  stronger version of (A1).\\

\begin{enumerate}
\item[(A1')] $|r_i(s,a)|\le w(s)$ for each $(s,a)\in \mathbb{K}$ and for all $i\in I_0$.
\end{enumerate}

\begin{enumerate}
\item[(A2)]
\begin{equation}
\label{a12}
\lim\limits_{n\to\infty} \sup_{\pi\in\Pi} E^\pi_\mu\left( \sum_{k=n}^\infty
\beta^{k -1}w(s_{k})\right)=0,
\end{equation}
and for each $n\in \mathbb{N}$
\begin{equation}
\label{a21}
\lim\limits_{l\to\infty}\sup_{\pi\in \Pi} E^\pi_\mu\left(w(s_n)1_{\{ w(s_n)\ge l\}}\right)=0
\end{equation}
\end{enumerate}

Assumption (\ref{a12}) is used in \cite{dufour} and is related to
condition (C) in \cite{schal1, fr}. Condition in (\ref{a21}) is a sort of the uniform integrability 
of the function $w$ with respect to strategic measures. 
It is related with Assumption {\bf B}  in \cite{dg}.
Clearly, here $1_{\{ w(s_n)\ge l\}}$ 
is understood as $1_{\{ w(s_n)\ge l\}}(s_n)$.  
\\

\noindent{\bf Remark 1} Observe that
\begin{eqnarray*}
0 &\le&
 \sup_{\pi\in\Pi}E_\mu ^\pi\left(\sum_{k=1}^\infty \beta^{k-1} w(s_k)1_{\{ w(s_k)\ge l\}}\right) \\ &\le&
\sup_{\pi\in\Pi}E_\mu ^\pi\left(\sum_{k=1}^{n-1} \beta^{k-1} w(s_k)1_{\{ w(s_k)\ge l\}}\right) +
 \sup_{\pi\in\Pi}E_\mu ^\pi\left(\sum_{k=n}^\infty \beta^{k-1} w(s_k)\right),
\end{eqnarray*}
for each $\pi\in\Pi.$  Take any $\epsilon >0$ and by (\ref{a12}) choose some $n$ such that
$$
0 \le
 \sup_{\pi\in\Pi}E_\mu ^\pi\left(\sum_{k=1}^\infty \beta^{k-1} w(s_k)1_{\{ w(s_k)\ge l\}}\right)\le
\sup_{\pi\in\Pi}E_\mu ^\pi\left(\sum_{k=1}^{n-1} \beta^{k-1} w(s_k)1_{\{ w(s_k)\ge l\}}\right) +\epsilon.
$$
This fact and (\ref{a21}) imply that
$$
0\le \limsup\limits_{l\to\infty} \sup_{\pi\in\Pi}E_\mu^\pi\left(\sum_{k=1}^\infty \beta^{k-1} w(s_k)1_{\{ w(s_k)\ge l\}}\right) \le \epsilon.
$$
Consequently, we have that
\begin{equation}
\label{a11}
\lim\limits_{l\to\infty}\sup_{\pi\in\Pi}E_\mu^\pi\left(\sum_{k=1}^\infty \beta^{k-1} w(s_k)1_{\{ w(s_k)\ge l\}}\right) =0.
\end{equation}
From (\ref{a21}) with $n=1$  or (\ref{a11}), it follows that   $\int_Sw(s)\mu(ds)<\infty.$

If the functions $r_i$ are bounded from above, then the function $w$ in (A1) may be constant.
In this case, (A2) is trivially satisfied. Some additional comments and examples,
where (A1) and (A2) are satisfied, are given in Section 6 and 7.\\

Under assumptions (A1)-(A2) all expectations in (\ref{CP}) are   well-defined, since
for  $r_i^+(s,a):=\max\{r_i(s,a),0\},$ $(s,a)\in \mathbb{K},$
by (\ref{a12}) we have
\begin{equation}\label{eqineqEF1}
0 \le J_i^+(\pi):= E_{\mu}^\pi
\left(\sum_{n=1}^\infty \beta^{n-1}r_i^+(s_n,a_n)\right) \le E_{\mu}^\pi\left(\sum_{n=1}^\infty \beta^{n-1}w(s_n)\right)<\infty.
\end{equation}
Let  $r_i^-(s,a):=\min\{r_i(s,a),0\},$ $(s,a)\in \mathbb{K},$  and
$$
J_i^-(\pi):= E_{\mu}^\pi
\left(\sum_{n=1}^\infty \beta^{n-1}r_i^-(s_n,a_n)\right).
$$

Therefore, every $J_i(\pi),$ $i\in I_0,$  can be written as
\begin{equation}\label{eqeqEF2}
J_i(\pi)= J_i^-(\pi)+ J_i^+(\pi).
\end{equation}

\section{Optimality of randomized stationary policies}

In order to formulate our main results we  need the following standard compactness and semicontinuity assumptions.
For MDPs with {\it weakly continuous transitions} we assume:\\

\begin{enumerate}
\item[(W1)] For each $s\in S,$ the action set $A(s)$ is compact and the set-valued mapping $s\to A(s)$
is upper semicontinuous, that is, the set $\{s\in S:\ A(s)\cap B\not= \emptyset\}$ is closed for every
closed set  $B\subset A.$
		\item[(W2)] For each bounded continuous function $v:S\to\mathbb{R}$, the function
$(s,a) \to \int_S v(t)p(dt|s,a)$ is continuous on $\mathbb{K}$ (weak continuity of $p$).
\item[(W3)]  The functions
  $r_i,$ $i\in I_0,$ are upper semicontinuous on $\mathbb{K}.$
\end{enumerate}

For MDPs with {\it setwise continuous transitions} we assume:

\begin{enumerate}
\item[(S1)] For each $s\in S,$ the action set $A(s)$ is compact
and the set $\mathbb{K}$ is Borel.
	\item[(S2)] For every $s\in S$ and for each $D \in {\cal B}(S)$,  the function
$p(D|s,\cdot)$ is continuous on $A(s)$ (setwise continuity of $p$).
\item[(S3)]  The functions
  $r_i(s,\cdot),$ $i\in I_0,$ are upper semicontinuous on $A(s)$ for each $s\in S.$
\end{enumerate}

We can now state our main result in this section.\\

\noindent {\bf Theorem 1} {\it Assume that } (A1)-(A2) {\it and either} (W1)-(W3) {\it  or}  (S1)-(S3)
{\it hold.  Then, there exists an optimal stationary policy $\varphi^*\in \Phi.$ }\\

Consider the following additional assumption.  
\begin{enumerate}
\item[(A3)] There exists $\delta>0$ such that
$$\int_Sw(t)p(dt|s,a)\le \delta w(s)\quad\mbox{for all}\quad (s,a)\in\mathbb{K}.$$
\end{enumerate}
If $\beta\delta <1,$ then (A3)  implies that (\ref{a12}) holds.\\

\noindent{\bf Corollary 1} {\it
Assume that } (A1) and (\ref{a21})  {\it and either} (W1)-(W3) {\it  or}  (S1)-(S3) {\it hold.  Assume in addition
that} (A3) is satisfied {\it and $\beta\delta<1.$
Then,  there exists an optimal stationary policy.} \\

\noindent{\bf Remark 2}
In models of MDPs (constrained or unconstrained) with weakly continuous transitions the following assumption is often made.
\begin{enumerate}
\item[(W4)]
The function $w$ in (A1) is continuous on $S$ and
$(s,a)\to\int_S w(t)p(dt|s,a)$ is continuous on $\mathbb{K}.$
\end{enumerate}
In models with setwise continuous transitions the following additional assumption is sometimes imposed.
\begin{enumerate}
 \item[(S4)]
For every $s\in S,$ the function
$a\to\int_S w(t)p(dt|s,a)$ is continuous on $A(s).$
\end{enumerate}
It should be noted that in general deterministic
policies are not sufficient for solving the problem (CP), see  for example \cite{frid}.\\

\noindent{\bf Remark 3}
Assumptions (A1) and (A3) were frequently used  in the studies
of discounted unconstrained MDPs with unbounded reward functions. The common approach is based on consideration of  the
{\it weighted norm} defined with the aid of the function $w$,   see \cite{w,vdw,hll}.
Under assumptions of Corollary 1, one can consider a standard transformation of the
CMDP  to a model with rewards bounded from above, see \cite{vdw}, Remark 2.5 in \cite{dufour}, or Section 10 in \cite{fp18}.
Making use of this transformation,  Corollary 1 can be deduced from Theorem 9.1 of \cite{fr},
 who studied CMDPs with reward functions $r_i$ bounded from above. However, one has to assume in addition
that (W4) or (S4) holds. That is because the transition probability functions in the transformed model depend on
$w$ and should have the weak or setwise continuity property.\\

\noindent{\bf Remark 4} The two alternative sets of conditions (W1)-(W3) and (S1)-(S3) were exploited in
the literature on   stochastic dynamic programming, see for example  \cite{b,bs,dy,hl,schal1,schal2}.
We would like to point out that there are several works on CMDPs with  transition probabilities
that are weakly continuous on $\mathbb{K}$ (jointly continuous in $(s,a)$)
and with  reward/cost functions that are upper/lower semicontinuous on $\mathbb{K}.$
Such assumptions are naturally satisfied in a number of models, where the  transition
law is determined by some difference equation with random shocks.
For a survey of various results on CMDPs with
weakly continuous transition probabilities   the reader is referred to
\cite{pbook,mp,gmmor,dufour,z} and \cite{fp1,fp2}.
CMDPs with both weakly and setwise continuous transitions
were  studied  in  \cite{fr}, whereas the models with measurable transition
and reward functions are examined in \cite{fp2,fp18}.

\subsection{The set of strategic measures}

Let $\Omega := (S\times A)^\infty.$ It is known that $\Omega$ is a Borel space. 
Let $\pi\in \Pi.$ We refer to $P_{\mu}^{\pi}$ as the {\it strategic probability measure}
generated by the policy $\pi$ and the initial distribution $\mu.$ Let
${\cal P}$ be the {\it set of strategic probability measures}, i.e.,
$${\cal P}:=\{P_{\mu}^{\pi}:\ \pi\in\Pi\}\subset \Pr(\Omega).$$
Under condition (S1), $\mathbb{K}^\infty$ is a Borel subset of $\Omega.$
Condition (W1) implies that $\mathbb{K}^\infty$ is a closed subset of $\Omega.$ 
For each $P^\pi_\mu \in {\cal P},$ $P^\pi_\mu(\mathbb{K}^\infty)=1.$
Let ${\cal C}_n$ (${\cal U}_n$) be the set of all bounded (bounded from above)
Borel functions on $(S\times A)^n $ having the following property.
A function $u$ belongs to ${\cal C}_n$ (${\cal U}_n$),
if  $u(s_1,\cdot,...,s_n,\cdot)$ is continuous
(upper semicontinuous) on $A(s_1)\times\cdots\times A(s_n)$ for any sequence of states $(s_1,...,s_n).$
Note that ${\cal U}_n$ contains the class of all upper semicontinuous and bounded from above functions
on $(S\times A)^n.$  \cite{schal1} defined the $ws^\infty$-topology on $\Pr(\Omega)$
as the coarsest topology in which
the functionals $P\to \int_\Omega u dP$ are continuous for each $u\in {\cal C}_n$ and $n\in \mathbb{N}.$
Assuming (S2), $A(s)=A$ for all $s\in S$ and $A$ is compact, \cite{schal1}  showed that $\cal P$ is compact
in the $ws^\infty$-topology on $\Pr(\Omega).$
\cite{schal2} also discussed a more general case, in which $A(s)$ may depend on $s.$ However, no formal proofs were given.
Later, \cite{now} observed that the relative $ws^\infty$-topology
on $\cal P$ is equivalent to the weak topology on this space of measures.\footnote{This result implies that the topology
on ${\cal P}$ is metrizable and one can think of convergence of sequences in ${\cal P}$ instead of nets.}
As a consequence,  the functionals
$P\to \int_\Omega udP$ are continuous on the compact space $\cal P$ endowed with the weak topology for any
$u\in {\cal C}_n$, $n\in \mathbb{N}.$  It is worthy to mention that assumption (S2),  
the initial distribution  and an argument related to the Scorza-Dragoni theorem (see \cite{ku}) 
play a fundamental  role  in the proof of \cite{now}.
\cite{b} extended the result of \cite{now} by allowing the action spaces
to depend on the partial histories of the process.

Let $\hat{{\cal U}}_n$ be the set  of all functions  from
${\cal U}_n$  restricted to the set $\mathbb{K}^n.$  Every
$u\in {\cal U}_n$  can be viewed as a function on $\mathbb{K}^\infty$ that 
depends only on the first $2n$ coordinates.

From Theorem 2.1 and Proposition 3.2 in \cite{b}, we conclude the following statement.\\

\noindent{\bf Lemma 1}
($a$) {\it Let assumptions }
(S1)-(S2) {\it or} (W1)-(W2) {\it be satisfied. Then, the set $\cal P$ is a compact subset of $\Pr(\Omega)$ endowed  with the weak topology.}\\
($b$) {\it If} (S1)-(S2) {\it hold, then the functional $P\to \int_{\mathbb{K}^\infty} udP$ is upper semicontinuous on $\cal P$
for each $u\in \hat{{\cal U}}_n,$ $n\in\mathbb{N}.$ } \\
($c$) {\it If} (W1)-(W2) {\it are satisfied, then  the functional $P\to \int_{\mathbb{K}^\infty} udP$ is  upper semicontinuous  on $\cal P$
for every  bounded  from above and  upper semicontinuous   function $u$ on $\mathbb{K}^n, $  $n\in\mathbb{N}.$ }\\

Note that equivalently the discounted reward functional  may be written as follows
\begin{equation}
\label{r}
 J_i(\pi)=  \int_{\mathbb{K}^\infty} \sum_{n=1}^\infty  \beta^{n-1} r^+_i(s_n,a_n)dP_{\mu}^\pi
+\int_{\mathbb{K}^\infty}\sum_{n=1}^\infty  \beta^{n-1} r^-_i(s_n,a_n)dP_{\mu}^\pi.
\end{equation}
Therefore, $J_i(\pi)$ can be viewed as a function of $P_{\mu}^\pi \in \cal P.$
Sometimes, we shall write $J_i(P^\pi_{\mu})$ for $J_i(\pi).$ \\
%%%%%%%%%%%%%%%%%%%%%%%%%%%%%%%%%%%%%%%%%%%%%%%%%%%%%%%%%%%%%%%%%%%%%%

 We now state our basic lemma.\\

\noindent{\bf Lemma 2}  {\it  Under assumptions of Theorem 1,
the discounted   reward functionals $J_i:{\cal P}\to \mathbb{R}_-$ 
are upper semicontinuous for all  $i\in I_0.$}\\

\noindent {\bf Proof }
First, we prove the lemma   assuming  (S1)-(S3) and that every function $r_i$ is non-negative.
Consider the truncated functions $r_i^l(s,a)=\min\{l,r_i(s,a)\}$ 
for $(s,a)\in \mathbb{K},$  $l\in\mathbb{N}$ and $i\in I_0.$
Then, every function $r_i^l(s,\cdot)$ is upper semicontinuous on $A(s)$ for every $s\in S.$
Note that  under condition (A1), for every $N\in \mathbb{N}$,
\begin{equation}\label{szog}
\sup_{\pi\in\Pi}
E_{\mu}^\pi\left(\sum_{n=N+1}^\infty \beta^{n-1} r_i^{l}(s_n,a_n)\right)\le
\sup_{\pi\in\Pi}
E_{\mu}^\pi\left(\sum_{n=N+1}^\infty \beta^{n-1} w(s_n)\right).
\end{equation}
Let
$$
J_i^l(\pi):= E_{\mu}^\pi\left(\sum_{n=1}^\infty \beta^{n-1} r_i^{l}(s_n,a_n)\right)\quad \mbox{and}\quad
J_i^{l,N}(\pi):= E_{\mu}^\pi\left(\sum_{n=1}^N \beta^{n-1} r_i^{l}(s_n,a_n)\right).$$
From (\ref{szog}) and (\ref{a12}), it follows that $J_i^{l,N}(\pi)$ converges to $J_i^l(\pi)$
as $N\to\infty,$ uniformly in $\pi\in\Pi.$ Therefore, it is sufficient to prove that for each
$N,$  $J_i^{l,N}(\pi)$ (understood as a function of $P_{\mu}^\pi$)
is upper semicontinuous on $\cal P.$
Observe that
\begin{eqnarray} \label{lhs} \nonumber
\lefteqn{
\sup_{\pi\in\Pi}\left| E_{\mu}^\pi\left(\sum_{n=1}^N \beta^{n-1} r_i(s_n,a_n)\right)-
E_{\mu}^\pi\left(\sum_{n=1}^N \beta^{n-1} r_i^{l}(s_n,a_n)\right)\right| \le} \\
&& \sum_{n=1}^N \sup_{\pi\in\Pi}  E_{\mu}^\pi \left(r_i(s_n,a_n)- r_i^{l}(s_n,a_n)\right)\le
\sum_{n=1}^N \sup_{\pi\in\Pi}  E_{\mu}^\pi \left( w(s_n)1_{\{w(s_n)\ge l\}}\right).
\end{eqnarray}
By (\ref{a21}) in (A2), for each $n=1,...,N,$
$$ 0\le  \sup_{\pi\in\Pi}  E_{\mu}^\pi\left(w^l(s_n)
1_{\{w(s_n)\ge l\}}\right)\to 0\quad \mbox{as}\quad l\to\infty.
$$
This fact and (\ref{lhs}) imply that
\begin{equation}
\label{limit0}
\sup_{\pi\in\Pi}\left| E_{\mu}^\pi\left(\sum_{n=1}^N \beta^{n-1} r_i(s_n,a_n)\right)-
E_{\mu}^\pi\left(\sum_{n=1}^N \beta^{n-1} r_i^{l}(s_n,a_n)\right)\right| \to 0 \quad\mbox{as}\quad l\to\infty.
\end{equation}

By Lemma 1$(b)$, the functional $P \to \int_{\mathbb{K}^\infty} u^l_i dP$ with
$$
u^l_i(s_1,a_1,\ldots,s_N,a_N)=\sum_{n=1}^N \beta^{n-1} r_i^{l}(s_n,a_n)
$$
is upper semicontinuous on $\cal P.$  Since the uniform limit of a sequence of 
upper semicontinuous functionals is upper semicontinuous,
using (\ref{limit0}), one can easily conclude that
$P \to \int_{\mathbb{K}^\infty} u_i dP$   with
$$u_i(s_1,a_1,\ldots,s_N,a_N)=\sum_{n=1}^N \beta^{n-1} r_i(s_n,a_n)$$
is also upper semicontinuous on $\cal P.$
Thus, we have proved the lemma under conditions (S1)-(S3)  for non-negative functions $r_i$, $i\in I_0.$
The proof for $r_i\ge 0$ under assumptions (W1)-(W3) makes use of  Lemma 1$(c)$ and proceeds along the same lines.  
Hence, we can conclude, in both cases under consideration, that the functionals
$P^\pi_\mu \to J_i^+(\pi)$ are upper semicontinuous on $\cal P.$
Since  $r^-_i \le 0,$ every functional
$P^\pi_\mu \to J_i^-(\pi)$ is also upper semicontinuous on $\cal P.$  The assertion now follows because
 $J_i (\pi)=J_i^-(\pi)+J_i^+(\pi).$
 $\Box$\\

In Section 6  we show a relation of Lemma 2 under conditions (W1)-(W3) with some recent results on weak convergence
of measures and unbounded mappings.

\subsection{ Occupancy measures and randomized stationary optimal policies}

Let $Q_\mu^{\pi}(ds\times da)$ be the {\it occupancy measure} on ${\cal B}(S\times A)$ of a policy $\pi\in\Pi,$
i.e., the measure defined as %follows
\begin{equation}
\label{om1}
Q_\mu^\pi(Z):= E_{\mu}^\pi\left(\sum_{n=1}^\infty\beta^{n-1}1_{Z}(s_n,a_n)\right),
\quad Z\in{\cal B}(S\times A).
\end{equation}
This  measure is finite and concentrated on the set  $\mathbb{K}$ for $\beta\in (0,1).$ Since every $\pi \in \Pi$ determines uniquely
the probability measure $P^\pi_{\mu} \in {\cal P},$ equality (\ref{om1}) shows how
$Q_\mu^\pi$ is determined by $P^\pi_{\mu}.$

Let ${\cal Q}_\mu^{\Pi}$ and ${\cal Q}_\mu^\Phi$
be the sets of all occupancy measures of policies $\pi\in \Pi$ and $\varphi\in\Phi,$ respectively.
By Lemma 4.1 in \cite{fr} or Proposition D8 in \cite{hl}, for any $\pi\in \Pi$ 
there exists some $\varphi\in \Phi$ such that
$$
Q_\mu^\pi(B\times D)= \int_B \varphi(D|s)q_\mu^\pi(ds),
\quad B\in{\cal B}(S),\ D\in{\cal B}(A),
$$
where $q_\mu^\pi$ is the projection of $Q_\mu^\pi$ on $S,$ i.e.,
$q_{\mu}^\pi(B)=Q_\mu^\pi(B\times A)$ for any $B \in {\cal B}(S).$
Moreover,  $Q_\mu^\pi =Q_\mu^\varphi.$ The proof of this fact is given  
in \cite{bor}  (see Lemma 3.1). Although \cite{bor} considered
models on a countable state space, his proof also applies to our framework, since it does not require
any continuity assumptions of the transition probability. The same result for Borel state space models
was reported in Lemma 24 in \cite{pbook} and Theorem 1 in \cite{z} or Lemma 4.2 in \cite{fr}.
Therefore, we can formulate the following  result
(see Corollary 4.1 in \cite{fr}).\\

\noindent{\bf Lemma 3}    ${\cal Q}_\mu^{\Pi}={\cal Q}_\mu^\Phi.$ \\

Lemma 3 directly implies  the following statement, on which the convex 
analytic approach to MDPs is based; see Remark 6 and, e.g.,
\cite{bor,pbook,mp,gmmor}.\\

\noindent{\bf Lemma 4}  {\it  For
 each $\pi\in \Pi$ there exists some $\varphi \in \Phi$ such that $J_i(\pi)=J_i(\varphi)$ for all $i\in I_0.$}\\

\noindent {\bf Proof of Theorem 1}  Lemmas 1 and  2  imply that
$${\cal P}^*:=\{P^\pi_{\mu}\in {\cal P}: \ J_1(P^\pi_{\mu})\ge d_1,\ldots, J_m(P^\pi_{\mu})\ge d_m \}$$
is a compact subset of ${\cal P}.$
Therefore, there exists a strategic measure $P^{\pi^*}_{\mu}\in {\cal P}^*$ such that
$$ \max_{P^{\pi}_{\mu}\in{\cal P}^*}J_0(P^{\pi}_{\mu})=J_0(P^{\pi^*}_{\mu}).$$
By Lemma 4,  there exists some
$\varphi^*\in \Phi$ such that $J_i(\pi^*)=J_i(\varphi^*)$
for all $i\in I_0.$ Clearly, $\varphi^*$ is a solution to problem (CP).
$\Box$\\

\noindent{\bf Remark 5}
Assume that $S$ contains an absorbing state $s^*$ with zero rewards. Then,
$p(S\setminus \{s^*\}|s,a)\le 1$ for all $(s,a)\in\mathbb{K}.$ Assumptions (A1) and (A2) can
be considered with $\beta =1$ and $w$ such that $w(s^*)=0$ and $w(s)\ge 1$ for all $s\not= s^*.$
 Lemma  2 with this modification remains correct.
If the other  assumptions of Theorem 1 are satisfied, then Lemmas 1 and 2
imply the existence  of an optimal policy. In general, Lemma 4 may not hold for $\beta=1;$
see examples in \cite{fso}. However, this lemma and the suggested  version of Theorem 1
with $\beta=1$  hold  for absorbing MDPs considered in  \cite{alt, fr}. \\

\noindent{\bf Remark 6} From (\ref{om1}), it follows that  
\begin{equation}
\label{om2}
\int_\mathbb{K} r(s,a)Q_\mu^\pi(ds\times da)= E_{\mu}^\pi\left(\sum_{n=1}^\infty\beta^{n-1}r(s_n,a_n)\right),
\end{equation}
for all non-negative measurable function $r:S\times A \to \mathbb{R}.$
Formulae (\ref{eqineqEF1}) and (\ref{eqeqEF2}) imply that (\ref{om2}) holds for $r=r_i,$ $i\in I_0.$
By Lemma 1$(a)$ and (\ref{om2}), we conclude  that ${\cal Q}_\mu^\Pi$ is compact
in the weak topology on the space of measures  $\nu$  on  ${\cal B}(S\times A)$
such that $\nu(S\times A)=1/(1-\beta).$
For a more detailed discussion the reader is referred to Section 4 in \cite{fr}.
Convexity of ${\cal Q}_\mu^\Pi$ follows from convexity of $\cal P$ using standard arguments based on disintegration of measures
on the product space, see, e.g., \cite{pbook}, \cite{schal2} or Corollary 4.1 in \cite{fr}.
We have shown that problem (CP) has a solution $\varphi^*\in \Phi.$
It is now clear that $Q_\mu^{\varphi^*}$ solves the following linear programming problem:
$$\mbox{(CP0)}\quad \mbox{maximize}\quad
\int_\mathbb{K} r_0(s,a)Q_\mu^\pi(ds\times da)$$
$$ \mbox{subject to}\quad
Q_\mu^\pi\in {\cal Q}_\mu^\Pi\quad \mbox{and}\quad
\int_\mathbb{K} r_i(s,a)Q_\mu^\pi(ds\times da)\ge d_i, \quad i\in I.
$$\\

\noindent{\bf Remark 7} Under  additional assumptions one can give a characterization of the solution
$Q_\mu^{\varphi^*}$ to problem (CP0) as in \cite{pbook,mp} or \cite{z}. The first  assumption is (A1').  The second condition requires that
 $r_i(s,\cdot)$ is continuous on  $A(s)$ for all $s\in S$, $i\in I.$ By Lemma 2 applied to $r_i$ and $-r_i,$  the functional
$P^\pi_\mu \to J_i(\pi)$  is continuous on $\cal P.$   
The last assumption is   Slater's condition demanding that there exists a policy $\pi' \in \Pi$ such that
$J_i(\pi') > d_i$ for all $i\in I.$ Using the Lagrange functional approach as in \cite{pbook} and \cite{z}, with minor
modifications, one can prove that $Q_\mu^{\varphi^*} = \sum_{j=1}^{m+1}\xi_jQ_\mu^{f_j},$
where $f_j\in F,$ $\xi_j\ge 0$ for all $j$ and $\sum_{j=1}^{m+1}\xi_j=1.$
Three ingredients play a significant role in the proof:  the 
relation  between $P^\pi_{\mu}\in {\cal P}$  and $Q_\mu^\pi$
with $\pi=\varphi$ given in (\ref{om1}), the   compactness and convexity of $\cal P$ and  finally, 
Lemma 3, which implies that $Q_\mu^\varphi \to \int_\mathbb{K} r_i(s,a)Q_\mu^\varphi(ds\times da)$
is continuous (upper semicontinuous) on ${\cal Q}_\mu^\Phi$ for each $i\in I$ (for $i=0$).
Linear programming problems for CMDPs with weakly continuous transitions and unbounded
cost functions satisfying condition similar to (A1) and some additional assumptions
were studied by \cite{dufour}. \cite{fr} and \cite{fp18}, on the other hand,
extended many earlier results on CMDPs to a class of
total reward models with absorbing states including discounted ones.

\section{Chattering and deterministic  optimal stationary policies}

Let $f_1,...,f_N \in F$ be any stationary deterministic policies
and $\alpha_1,...,\alpha_N$ be non-negative numbers such that $\sum_{j=1}^N \alpha_j =1.$ By $Q^{f_j}_\mu$
we denote the occupancy measure
induced by $f_j\in F.$ Its projection on $S$ is denoted by $q^{f_j}_\mu.$ Define
$$Q(\cdot):= \sum_{j=1}^N \alpha_j Q^{f_j}_\mu(\cdot).$$
Note that $Q $ is also an occupancy measure since the set of occupancy measures is convex.   
By Proposition D8 in \cite{hl}, there exists
 $\varphi\in \Phi$ such that 
\begin{equation}
\label{wazne}
Q(B\times D)= \int_B \varphi(D|s)q(ds),
\quad\mbox{for all}\quad  B\in{\cal B}(S),\ D\in{\cal B}(A),
\end{equation}
where $q$ is the projection of $Q$ on $S.$
Since $Q$ is an occupancy measure, we now write  
$Q=Q^\varphi_\nu$ and $q=q^\varphi_\mu$ in (\ref{wazne}).
Note that for any bounded Borel measurable function $g:\mathbb{K}\to
\mathbb{R}$   
$$\int_{\mathbb{K}}g(s,a)Q^\varphi_\mu(ds\times da)= \int_S\int_{A(s)}g(s,a)\varphi(da|s)q^\varphi_\mu(ds).$$
If $\delta_{f_j(s)}(\cdot)$ is the Dirac measure with support at the point
$f_j(s),$ then
$$\int_{\mathbb{K}}g(s,a)Q^{f_j}_\mu(ds\times da)= \int_S\int_{A(s)}g(s,a)\delta_{f_j(s)}(da)q^{f_j}_\mu(ds)
= \int_Sg(s,f_j(s))q_\mu^{f_j}(ds).$$
\\

\noindent{\bf Lemma 5} {\it There exist Borel measurable functions $\gamma_j:S \to [0,1]$  ($j=1,...,N$) such that
 $\sum_{j=1}^N\gamma_j(s)=1$ and $\varphi (\cdot|s)= \sum_{j=1}^N \gamma_j(s)\delta_{f_j(s)}(\cdot)$
\ ($q_\mu^\varphi$-a.e.}). \\

\noindent{\bf Proof}
Observe that $q^\varphi_\mu = \sum_{j=1}^N \alpha_j q^{f_j}_\mu.$  Hence, 
$q^{f_j}_\mu \ll q^\varphi_\mu$ for  every $j=1,\ldots,N.$  Let $\rho_j:=\frac{q^{f_j}_\mu}{q^\varphi_\mu}$
be a non-negative Borel measurable Radon-Nikod{\'y}m derivative of
$q^{f_j}_\mu$ with respect to $q^\varphi_\mu.$
Define $\gamma_j(s)=\alpha_j\rho_j(s),$ $s\in S.$
Observe that for every  $B\in{\cal B}(S)$  
$$\int_B\left(\sum_{j=1}^N\gamma_j(s)\right) q^\varphi_\mu(ds) =
\int_B\left(\sum_{j=1}^N \alpha_j\rho_j(s)\right)q^\varphi_\mu(ds)= \sum_{j=1}^N\alpha_j
q^{f_j}_\mu(B)= q^\varphi_\mu(B).$$
This implies that $\sum_{j=1}^N\gamma_j(s)=1$
for all $s\in B_1$ where $B_1 \in {\cal B}(S)$ and $ q^\varphi_\mu(B_1)=1.$ For any $s\notin B_1$ we can
modify our definition of $\gamma_j(s).$   Namely, we can put
$\gamma_1(s)=1$ and $\gamma_j(s)=0$ for $j=2,...,N.$

Let  
\begin{equation}
\label{varp}
\overline{\varphi}(\cdot|s):= \sum_{j=1}^N \gamma_j(s)\delta_{f_j(s)}(\cdot).
\end{equation}
For a bounded Borel measurable function $g$  on $\mathbb{K},$  we have
\begin{eqnarray}
\label{wazne3}
\nonumber
\lefteqn{\int_S\int_{A(s)} g(s,a)\varphi(da|s)q^\varphi_\mu(ds)=\int_{\mathbb{K}}g(s,a)Q^\varphi_\mu(ds\times da)} \\  &=&
\sum_{j=1}^N\alpha_j\int_{\mathbb{K}}g(s,a)Q^{f_j}_\mu(ds\times da) \nonumber
=\sum_{j=1}^N\int_S\int_{A(s)}\alpha_jg(s,a)\delta_{f_j(s)}(da)q^{f_j}_\mu(ds) \\  & =&
\int_S\int_{A(s)}g(s,a)\overline{\varphi}(da|s)q^\varphi_\mu(ds).
\end{eqnarray}
Since (\ref{wazne3}) holds for every bounded Borel measurable function $g:\mathbb{K}\to\mathbb{R},$
we conclude that $\varphi(\cdot|s)=\overline{\varphi}(\cdot|s)$ \ ($q^\varphi_\mu$-a.e.).  $\Box$\\

The following terminology is borrowed from the theory of variational
calculus and control theory, see \cite{ro}.
A stationary policy $\phi \in \Phi$ is called {\it chattering}, if there exist a family of $N$ Borel functions
$\gamma_j:S\to [0,1]$  and a family of $N$ deterministic stationary policies  $f_j \in F$ such that
$$\phi(D|s)= \sum_{j=1}^N\gamma_j(s)\delta_{f_j(s)}(D) \quad\mbox{and}\quad \sum_{j=1}^N\gamma_j(s)=1$$
for each $D \in {\cal B}(S)$ and  for all
$s\in S.$

Let $\pi\in\Pi.$  Following \cite{fr}, we define the {\it performance vector} ${\cal  V}(\pi):= (J_0(\pi),J_1(\pi),...,J_m(\pi))$ and the
{\it performance set} ${\cal V}:= \{ {\cal V}(\pi):  \pi\in\Pi\}.$\\

We are now ready to state our main result in this section.\\

\noindent{\bf Theorem 2} {\it Under assumptions of Theorem 1,
there exists a chattering stationary  policy} $\phi \in \Phi$ {\it  with} $N=m+1$
{\it that solves problem } (CP).\\

\noindent
{\bf Proof} By Theorem 1, there exists an optimal stationary policy $f^*\in \Phi.$
Using the same geometric arguments as in  the proof of Theorem 9.2 in \cite{fr}, one can conclude that
${\cal V}(f^*)$ lies on the boundary of the performance set $\cal V$ and therefore  the occupancy measure $Q^{f^*}_\mu$
can be represented as
$$
Q^{f^*}_\mu = \sum_{j=1}^{m+1} \alpha_jQ^{f_j}_\mu
$$  with some $f_1,...,f_{m+1} \in F$
and   non-negative numbers  $\alpha_1,...,\alpha_{m+1}$ such that $\sum_{j=1}^{m+1} \alpha_j =1.$
Let $\varphi\in\Phi$ be  as in (\ref{wazne}) with $Q=  Q^{f^*}_\mu.$
By Lemma 5, there exist Borel measurable functions
$\gamma_j:S\to [0,1]$ ($j=1,...,m+1$) such that $\sum_{j=1}^{m+1}\gamma_j(s)=1$ for all
$s\in S$, and a family of $m+1$ deterministic stationary policies   $f_j \in F$ such that for
$$\phi(\cdot|s):= \sum_{j=1}^{m+1}\gamma_j(s)\delta_{f_j(s)}(\cdot)$$
we have $\phi(\cdot|s)=\varphi(\cdot|s)$ \ ($q^\varphi_\mu$-a.e.). 
Moreover, it follows that $Q^{f^*}_\mu= Q^\varphi_\mu= Q_\mu^\phi.$ This implies that
$J_i(f^*) =J_i(\phi)$ for all $i\in I_0.$   Since $f^*$ is optimal, 
the chattering stationary policy $\phi$ is optimal as well. $\Box$ \\

\noindent{\bf Remark 8} \cite{fr} showed in Theorem 9.2(i) that for any feasible  policy $\pi$
there exist $m+2$ stationary deterministic policies $f_j\in F$ and non-negative numbers $\alpha_1,...,\alpha_{m+2}$
such that
$$Q^\pi_\mu =\sum_{j=1}^{m+2} \alpha_j Q_\mu^{f_j}$$
and $\sum_{j=1}^{m+2}\alpha_j=1.$ Their proof can also be
used under assumptions of Theorem 1. 
Applying Lemma 5 with $N=m+2$ one can easily conclude that $\overline{\varphi}$ defined in (\ref{varp})
has the property that
$Q_\mu^\pi =Q^\varphi_\mu =Q_\mu^{\overline{\varphi}}.$
Therefore, for the chattering policy $\overline{\varphi}$ in (\ref{varp}) with $N=m+2$
 we have that $J_i(\pi)=J_i(\varphi) =J_i(\overline{\varphi})$ for all $i\in I_0.$ \\

Combining Corollary 1  with a general ``purification result'' stated as Corollary 10.2
in \cite{fp18}, we can conclude the following fact.\\

\noindent{\bf Corollary 2} {\it Assume that  the initial distribution and transition probabilities
are atomless and assumptions of Corollary 1
with} (A1) {\it replaced by} (A1') {\it are satisfied.
Then, there exists a deterministic
stationary  policy $\tilde{f}  \in F$ that solves problem } (CP).\\

\noindent{\bf Proof}
Under  assumptions of Corollary 1 with (A1') instead of (A1), the model can be transformed
to an absorbing CMDP with bounded reward functions. This is mentioned in Section 10 of \cite{fp18}.
By Corollary 1, there exists an optimal  stationary policy  $f^*$ in the original model. (Since no continuity
conditions are imposed on the function $w$, we cannot in this place conclude the existence of $f^*$
in the transformed model established, e.g., in \cite{fr}.)  By Corollary 10.2 in \cite{fp18},
there exists some deterministic stationary policy $\tilde{f}\in F$
giving the same expected discounted rewards in both transformed and original model. Thus, we have
$J_i(f^*)=J_i(\tilde{f})$ for all $i\in I_0.$ Obviously, $\tilde{f}$ is an optimal policy.
$\Box$

 \section{ Examples }

In this section, we provide three examples satisfying assumptions imposed in  Corollaries 1-2.
In addition, we indicate a class of examples in which for any $\beta<1$ there exists some $\delta >0$ such that (A3) is satisfied
and $\delta\beta<1.$ This fact implies that (\ref{a12}) is satisfied. By Lemma 9 in Section 7 assumption (A3)
and some continuity assumptions imply that (\ref{a21})  holds.

We start with two  models of economic growth
theory, see \cite{stach}. In Example 1 the dynamics is deterministic, whereas Example 2 includes a stochastic component.   \\

\noindent{\bf Example 1} Consider  a  dynamic growth model with $S=[0,+\infty)$  and $A(s)=[0,s].$ Assume that
the level of the resource stock evolves according to the equation   $s_{n+1}=\sqrt{s_n-a_n},$ where $n\in\mathbb{N}.$
Assume that the utility or reward function for the economic agent  is $r_0(s,a)=a-1/a$ for $(s,a)\in \mathbb{K},$ with $a\in A(s):= [0,s],$
whilst the reward function of the authorities is $r_1(s,a)=\ln s$ for $(s,a)\in \mathbb{K}$. For $s=a=0,$ $r_1(s,a):= -\infty.$   
The agent's problem is  to solve problem (CP)
with the  constant $d_1$  provided by the authorities.

Let $w(s)=s+c$ for $s\in S$ with some constant $c\ge 1$ and let $\mu$ be an   initial distribution on $S$  such that $\int_Sw(s)\mu(ds)<\infty.$
Observe that (A1) holds. Moreover,  it must hold
$$w(\sqrt{s-a})= \sqrt{s-a}+c\le \delta (s+c) \ \mbox{for all}\ a\in [0,s] \ \mbox{and }\ s\in S$$
with some $\delta>0.$ Note that, for all $(s,a)\in \mathbb{K},$
we have
$$\frac{\sqrt{s-a}+c}{s+c}\le\frac{\sqrt{s}+c}{s+c}= \le1+\frac{\sqrt{s}-s}{1+c}\le 1+\frac 1{4+4c}.$$
Put $\delta:=1+1/(4+4c).$ It is easily seen that for every $\beta<1$ there exists $c\ge 1$ so that $\beta\delta<1.$
Obviously, (W1)-(W4) are satisfied and by Theorem 2 there exists an optimal stationary chattering policy.\\

\noindent
Based on the aforementioned example, we may generalize the method for finding $\delta$ such that $\delta\beta<1$
for the given discount factor $\beta<1.$
Suppose that we have some continuous function $w_0:S\to[1,+\infty)$ such that it holds
$$\int_S w_0(t)q(dt|s,a)\le  w_1(s), \quad (s,a)\in \mathbb{K},$$
where  $w_1$ is a non-negative continuous function such that $\theta:=\sup_{s\in S} (w_1(s)-w_0(s))<\infty.$
Then, we  define $w(s):=w_0(s)+c$ for $s\in S$ and some $c\ge 1.$ Simple calculations give that
$$\frac{\int_S(c+ w_0(t))q(dt|s,a)}{w_0(s)+c}\le \frac{w_1(s)+c}{w_0(s)+c}\le1+\frac{w_1(s)-w_0(s)}{1+c}$$
for any $(s,a)\in\mathbb{K}.$
Define now
$$\delta:=\max\{1+\frac \theta{1+c},1\}.$$
Hence, if $\theta>0$ for every $\beta<1$ we may take $c$ sufficiently large in order to have $\beta\delta<1.$ In the second case,
when $\theta\le 0,$ the condition $\beta\delta<1$ is always satisfied.  Many examples, which we have encountered in the literature
can be reduced to this case, when
\begin{itemize}
\item $S=[0,\infty)$,
 \item the functions $w_0$ and $w_1$ are increasing,
 \item  equation   $w_0(s)=w_1(s)$ has a unique solution  $s^*>0.$
 \item $w_1(s)-w_0(s)<0$ for all $s>s^*.$
\end{itemize}

The next two examples are dedicated to the application of Corollary  2.
The first example is given with weakly continuous transition probabilities, whereas in
the second one the transition probabilities are setwise continuous.\\

\noindent{\bf Example 2} Consider the model from Example 1, but with different dynamics  for the resource stock.
Suppose that the level of the resource is
described by the following equation
$$s_{n+1}=s_n-a_n+\sqrt{s_n-a_n}+\xi_n, \quad\mbox{for } n\in\mathbb{N}.$$
where $(\xi_n)_{n\in\mathbb{N}}$ is a sequence of i.i.d. random variables taking values in the interval $[0,+\infty).$ 
Moreover,  every
$\xi_n$ has an atomless distribution $\rho$ such that $\overline{m}:=\int_0^\infty z\rho(dz)<\infty.$
Additionally, assume the initial state is chosen at random according to an atomless  measure $\mu$. 
Suppose that $r_0(s,a)=\sqrt{a}$ and   
$r_1(s,a)=\ln(s+1)$ for every
$(s,a)\in \mathbb{K}.$ Let $d_1$ be a given number. The agent again faces problem (CP).
Let us define $w(s)=s+c$ with some $c\ge 1.$ Then, (A1') is satisfied. Furthermore, we obtain
$$\int_S w(t) q(dt|s,a) = s-a+\sqrt{s-a}+c +\overline{m}\le s+\sqrt{s}+c+\overline{m},\quad\mbox{for } (s,a)\in \mathbb{K}.$$
Consequently,
$$\frac{\int_S w(t) q(dt|s,a)}{s+c}\le 1+\frac{\sqrt{s}
+\overline{m}}{s+c}\le 1+\frac12\frac1{\sqrt{\overline{m}^2+c}-\overline{m}}.$$
Thus, (A3) is satisfied with $\delta:= 1+\frac12\frac1{\sqrt{\overline{m}^2+c}-\overline{m}}.$ Note that
for any value of the discount coefficient $\beta<1,$ we may choose
$c\ge 1$    so that $\beta\delta<1.$ Hence, by Corollary 2, 
there exists an optimal deterministic  stationary policy solving (CP).\\

\noindent{\bf Example 3} 
Assume that a system can be in a state $s\in [0,1]$,
where $0$ denotes the perfect condition and $1$ means that the system is completely broken.  
At the end of each month the system is checked, and its state is observed.  The higher values of $s,$ 
the worse condition of the system is. The controller each month decides
about  the repair. He chooses some $a\in A:=[0,1],$ where $a=0$ means no repair and $a=1$ is the replacement of the old system
by a new one. Hence, the larger $a,$ the more serious repair is required.
The transition probability $p$ is absolutely continuous with respect to the Lebesgue measure on $S$. In other words,
there exists a density $g(\cdot,s,a)$ for every $(s,a)\in \mathbb{K}$ such that
$$p(B|s,a)=\int_B g(t,s,a)dt\quad \mbox{ for } (s,a)\in\mathbb{K} \quad\mbox{and } B\in{\cal B}(S).$$
Additionally, let $g(t,s,\cdot)$ be continuous on $A$ for every $t,s\in S.$  The cost associated with the repair is
$c_0(s,a)$ for $(s,a)\in\mathbb{K}.$ Moreover, $c_0(s,\cdot)$ is bounded and lower semicontinuous on $A$ for every $s\in S.$
The management of the company requires that the sum of the discounted values describing the
system's state  cannot be  greater that $d_1$ (i.e., $c_1(s,a)=s$ for $(s,a)\in\mathbb{K}$).
The initial distribution $\mu$ is atomless.  Hence, the constrained control problem is
$$\mbox{(CP0)}\qquad \mbox{minimize}\quad E_\nu^\pi\left(\sum_{n=1}^\infty \beta^{n-1}c_0(s_n,a_n)\right)$$
 $$ \mbox{subject to}\quad
 E_\nu^\pi\left(\sum_{n=1}^\infty \beta^{n-1}s_n\right)\le d_1.
$$
Note that all assumptions in Corollary 2 are satisfied.
Thus, there exists an optimal deterministic stationary policy for problem (CP0).

\section{A comment on the approach of Dufour and Prieto-Rumeau}

In this section, we show that assumptions imposed by \cite{dufour} imply the uniform integrability condition (\ref{a21})
and that the $w^\alpha$-topology introduced in their paper is equivalent with the standard weak topology on $\cal P.$
We also demonstrate how Lemma 2 under assumptions
(A1)-(A2) and (W1)-(W3)  can be briefly deduced from recent results on weak convergence of finite measures.

Let $\cal M$ be a family of finite measures on a Borel space $Y.$
An extended real-valued Borel function
$v: Y\to  [-\infty,\infty]$ is called uniformly integrable  with respect to the family ${\cal M}$ if
\[
\lim_{l\to \infty} \sup_{\eta\in{\cal M}} \int_Y |v (y)|
		1_{ \{ |v (y)| \geq l\}} \eta (dy) = 0.
		%\label{eq:tv:ui}
\]
The following known lemma is useful for some considerations  in this paper. \\

\noindent{\bf Lemma 6 } {\it For a  sequence $(\eta_n)$ of finite measures on a Borel space $Y$
converging weakly to a measure $\eta$ and for
an upper  semicontinuous function $v:Y\to [-\infty,\infty],$  the following statements hold:\\
$(a)$ if $v^+=\max\{v,0\}$   is uniformly integrable with respect to
${\cal M}= \{\eta_n: \ n\in\mathbb{N}\},$ then $$\int_Y  v (y)\eta(dy)
\ge \limsup\limits_{n\to\infty}\int_Y v (y)\eta_n (dy),$$
$(b)$ if $v$   is continuous and
uniformly integrable with respect to
${\cal M}= \{\eta_n\},$  then the inequality in $(a)$ can be replaced by the equality. } \\

 Lemma 6$(a)$ follows from Fatou's lemma for uniformly integrable functions and weakly converging measures
 (Theorem 2.4 in \cite{fkl}). Statement $(b)$ is a particular case of Theorem 2 in \cite{za}. A related result was given in \cite{dg}, see
 Theorem 3.2.  Statement $(b)$ can be also concluded from point $(a)$ applied to the functions $v$ and $-v.$  \\

\cite{dufour} consider a CMDP that satisfies conditions (W1)-(W3), (A1), (\ref{a12}) (their Assumption  A, (B.3), C)
and some additional assumptions, called (B.1) and (B.2), involving a continuous function $w.$
They endow   the set of strategic measures $\cal P$ with the so-called $w^\alpha$-topology
($\alpha$ is a fixed discount factor). From their definition, it follows that
the $w^\alpha$-topology on $\cal P$  is  finer than the weak topology
and by  Theorem 3.10 in \cite{dufour}, $\cal P$ is a
metrizable compact space. From the   definition of the $w^\alpha$-topology on $\cal P,$
it follows that  the functional $P^\pi_\mu \to E^\pi_\mu(w(s_n))$ is continuous for any $n\in\mathbb{N}.$
 Since $w(s)1_{\{w(s) <l\}}(s)$ is bounded and lower semicontinuous on $S$,  the functional
$P^\pi_\mu \to E^\pi_\mu(w(s_n)1_{\{w(s_n) <l\}})$
is lower semicontinuous on $\cal P $  endowed with the weak topology.
Therefore,  it is lower semicontinuous on $\cal P$ endowed with the $w^\alpha$-topology. Thus, 
$$P^\pi_\mu \to E^\pi_\mu(w(s_n)1_{\{w(s_n) \ge l\}}) =
E^\pi_\mu(w(s_n)) -  E^\pi_\mu(w(s_n)1_{\{w(s_n) <l\}})$$
is upper semicontinuous on the compact  space
$\cal P $ in the $w^\alpha$-topology.
Note that
$$E^\pi_\mu(w(s_n)1_{\{w(s_n) \ge l\}}) \downarrow 0\quad\mbox{as} \quad l\to\infty.$$
From Dini's theorem  (see Chapter 9 in  \cite{royden}), it follows that
$$\lim_{l\to\infty}\sup_{\pi\in\Pi} E^\pi_\mu(w(s_n)1_{\{w(s_n) \ge l\}})=0,$$
that is, our assumption (\ref{a21}) is satisfied.
Theorem 1  shows that the standard weak topology on $\cal P$
is sufficient to solve problem (CP) under conditions
given by \cite{dufour}.

Let us assume that (\ref{a21}) holds and $w$ is continuous.
As follows from Lemma 6$(b)$, the functional $P^\pi_\mu  \to E^{\pi}_\mu(w(s_n))$ is continuous
on $\cal P$ endowed with the weak topology.
Hence, it follows that  the  weak topology is finer than  the $w^\alpha$-topology on $\cal P.$ Thus, the two topologies on $\cal P$ are equivalent.
%%%%%%%%%%%%%%%%%%%%%%%%%%%%%%%%%%%%%%%%%%%%%%%%%%%%

We now show how Lemma 2 under assumptions (A1)-(A2) and (W1)-(W3) can be deduced from Lemma 6.
Let us define for each $i\in I_0$  the function $u_i(s_1,a_,s_2,a_2,\ldots):=\sum_{n=1}^\infty\beta^{n-1} r_i(s_n,a_n)$
if $ \sum_{n=1}^\infty\beta^{n-1} r^+(s_n,a_n)<\infty,$ and
$u_i(s_1,a_,s_2,a_2,\ldots):=\infty$ otherwise.  This function is upper semicontinuous on
$\mathbb{K}^\infty.$  Formula (\ref{a11}) implies that $u_i$  are
uniformly integrable with respect to $\cal P.$
Lemma 6$(a)$ implies upper semicontinuity of the functional $P^\pi_\mu\to\int_{\mathbb{K}^\infty} udP_\mu^\pi.$

\section{Additional remarks on our basic assumptions}

Conditions (W1)-(W2), (W4) and (A3) imposed in Lemma 9 below, are discussed in Remark 2.3 in \cite{dufour}.
They imply their assumption (B.2) which plays an important role in proving compactness
of $\cal P$ in the $w^\alpha$-topology. This result was used in Section 6 in our proof that the assumptions of \cite{dufour}
imply the uniform integrability condition (\ref{a21}). Below we show directly that (W1)-(W2), (W4) and (A3) imply (\ref{a21})
(see the proof of Lemma 9). We also discuss the alternative case with conditions (S1)-(S2), (S4) and (A3).

Let us define two classes of non-negative  functions on $S$ denoted by $\hat{B}_+(S)$ and $\hat{U}_+(S)$, respectively.
A non-negative function $v$ belongs to the class $\hat{B}_+(S)$ ($\hat{U}_+(S)$) if it is Borel measurable
(upper semicontinuous) on $S$ and there exists a constant $c>0$ such that $v(s)\le c w(s)$  for all $s\in S.$
For any $v\in \hat{B}_+(S)$ define
\begin{equation}
\label{M}
(Mv)(s): =\sup_{a\in A(s)}\int_S v(t)p(dt|s,a), \  \ s\in S.
\end{equation}

\noindent{\bf Lemma 7} {\it Assume that } (A3) {\it is satisfied.}\\
$(a)$ {\it If $v\in \hat{U}_+(S)$ and } (W1)-(W2) {\it and} (W4) {\it hold, then $Mv \in \hat{U}_+(S).$}\\
$(b)$ {\it If $v\in \hat{B}_+(S)$ and } (S1)-(S2) {\it  and} (S4)  {\it hold, then $Mv \in \hat{B}_+(S).$}\\

\noindent{\bf Proof} $(a)$ There exists some  $c>0$ such that $cw-v\ge 0.$   The function $cw-v$ is
lower semicontinuous and non-negative. Thus,  it is a pointwise limit of a non-decreasing sequence of bounded  continuous functions.
Therefore,  from (W2), it follows  that  the function $(s,a)\to \int_S(cw(t)-v(t))p(dt|s,a)$ is lower semicontinuous
on $\mathbb{K}.$  This fact and (W4) imply that $(s,a)\to \int_Sv(t)p(dt|s,a)$  is upper semicontinuous on $\mathbb{K}.$
From (A3), it follows that there exists some constant $c_1 >0$ such that $(Mv)(s)\le c_1w(s)$ for all $s\in S.$
Using (W1) and the maximum theorem of \cite{berge},
we conclude that $Mv\in \hat{U}_+(S).$

$(b)$  The function $cw-v$ is non-negative for some $c>0$ and
is a pointwise limit of a non-decreasing sequence of bounded Borel measurable
functions. By (S2),   the function $a\to \int_S(cw(t)-v(t))p(dt|s,a)$ is lower semicontinuous
on $A(s)$ for each $s\in S.$  This fact and (S4) imply that $s\to \int_Sv(t)p(dt|s,a)$  
is upper semicontinuous on $A(s)$ for each $s\in S.$
By (A3),    $(Mv)(s)\le c_1w(s)$ for all $s\in S $ and for some $c_1>0.$  
Using Corollary 1 in  \cite{bp} and (S1), we conclude that $Mv\in \hat{B}_+(S).$
$\Box$ \\

\noindent{\bf Lemma 8} {\it Let $(v^l)$ be a non-increasing sequence
of functions on $S$ such that $v^l(t)\downarrow 0$ for each $t\in S$ as $l\to\infty.$
Let } (A3)  {\it be satisfied. Assume that  $v^l\in \hat{U}_+(S)$ for all $l\in \mathbb{N}$ and
conditions }  (W1)-(W2) {\it and } (W4) {\it hold, or  $v^l\in \hat{B}_+(S)$ for all $l\in\mathbb{N}$ and assumptions}
(S1)-(S2) {\it  and} (S4) {\it are satisfied.
Then, $(Mv^l)(s) \downarrow 0$
for each $s\in S$ as $l\to\infty.$}\\

\noindent{\bf Proof} From the proof of Lemma 7, we know that  in both cases, the function $a\to\int_Sv^l(t)p(dt|s,a)$
is upper semicontinuous on $A(s)$ for each $s\in S$, $l\in\mathbb{N}.$
By the monotone convergence theorem, $\int_Sv^l(t)p(dt|s,a) \downarrow 0$
for all $(s,a)\in \mathbb{K}$ as $l\to\infty.$ The assertion follows now from Dini's theorem. $\Box$\\

\noindent {\bf Lemma 9} {\it Let } (A3)  {\it be satisfied.
Assume that either
conditions }  (W1)-(W2) {\it and} (W4) {\it or} (S1)-(S2) {\it and }
(S4) {\it are satisfied. Then, condition }  (\ref{a21}) {\it in  } (A2) {\it holds.}\\

\noindent{\bf Proof}
Let $v^l(s_n)= w(s_n) 1_{\{w(s_n)\ge l\}}.$ For $n=1$ we have
$$0\le E_\mu^\pi(v_l(s_1))\le \int_S v^l(s_1)\mu (ds_1) \to 0
\quad \mbox{as} \quad l\to\infty.$$
If $n\ge 2$, then
$$ 0\le  \sup_{\pi\in\Pi}  E_{\mu}^\pi(v^l(s_n))
\le \int_S (M^{n-1}v^l)(s_1) \mu(ds_1),$$
where $M^{n-1}$ is the $(n-1)$st composition of the operator $M$ defined
in (\ref{M}) with itself.
By Lemmas 7 and  8 and the monotone convergence theorem, it follows by induction on $n$ that
$$ \lim_{l\to\infty}\int_S(M^{n-1}v^l)(s_1) \mu(ds_1) =0.$$
Thus, (\ref{a21}) in (A2) follows.  $\Box$\\

We now consider MDPs similar to those studied in \cite{jnjmaa}.
Let $(X_k)$ be a sequence of non-empty Borel subsets of S such that $X_k\subset X_{k+1}$ for each $k\in\mathbb{N}$
and $\bigcup_{k=1}^\infty X_k= S.$
Let $m_k= \sup_{s\in X_k}w(s),$ $k\in\mathbb{N}.$

\begin{enumerate}
\item[(A4)] For each $x\in X_k,$ $a\in A(s),$ $k\in\mathbb{N}$, $p(X_{k+1}|s,a)=1$  and $\sum\limits_{k=1}^\infty m_k\beta^{k-1}<\infty.$
\end{enumerate}

\noindent{\bf Lemma 10} {\it Under assumption} (A4) {\it conditions}
 (A1) {\it and} (A2) {\it hold with  $w(s):= m_k$ for all
$s\in X_k$, $k\in \mathbb{N}.$}\\

\noindent{\bf Proof} Clearly, (A1) and (\ref{a12})  are obviously satisfied. 
By the monotone convergence theorem, the condition in (\ref{a21}) holds for $n=1.$ Assume that $n\ge 2.$
For any non-negative Borel measurable function $\tilde{w}$ on $S$, we define by $E^\pi_{s_1}(\tilde{w}(s_n))$
the conditional expectation of $\tilde{w}$ with respect to the $n$th state given the initial state $s_1.$
Then, we have   $E^\pi_{\mu}(\tilde{w}(s_n))=\int_SE^\pi_{s_1}(\tilde{w}(s_n))\mu(ds_1).$
Choose any $s_1\in X_k,$ $k\in\mathbb{N}.$ Condition (A5) implies that $E^\pi_{s_1}(1_{X_{k+n-1}}(s_n))=1$ for all $\pi\in\Pi.$
Therefore, for each $l> m_{k+n-1}$ , we have
$\sup_{\pi\in\Pi}E^\pi_{s_1}(w(s_n)1_{\{w(s_n)\ge l\}})=0.$
Thus, $\lim_{l\to\infty}\sup_{\pi\in\Pi}E^\pi_{s_1}(w(s_n)1_{\{w(s_n)
\ge l\}})=0.$ By the monotone convergence theorem, we have
$\lim_{l\to\infty}\int_S\sup_{\pi\in\Pi}E^\pi_{s_1}(w(s_n)1_{\{w(s_n)\ge l\}})\mu(ds_1)=0.$
Hence,
$$
0\le \lim_{l\to\infty}\sup_{\pi\in\Pi}E^\pi_\mu(w(s_n)1_{\{w(s_n)\ge l\}})\le
\lim_{l\to\infty}\int_S\sup_{\pi\in\Pi}E^\pi_{s_1}(w(s_n)1_{\{w(s_n)\ge l\}})\mu(ds_1)=0,$$
which completes the proof of (\ref{a21}). $\Box$\\

Lemma 9 implies that (A2) holds, if  (A3) is satisfied and either $\beta\delta <1$ or  (A4) holds.
Lemma 10 provides an additional set of conditions under which (A1) and (A2) are satisfied.

Finally, we provide an example  illustrating the importance of our assumptions. In particular,
we show that if (\ref{a21}) does not hold, then an optimal policy may not exist even in a very simple model,
where (\ref{a12}) is trivially satisfied.\\

\noindent{\bf Example 4}  Consider the following unconstrained MDP, where
\begin{itemize}
\item $S=\{0,0^*,0^{**}\}\cup\mathbb{N}$;
\item $A(0^*)=\{0\}\cup\bigcup_{n=1}^\infty \{1/n\}$ and
 $A(s)=\{a_0\}$ for all $s\in S\setminus \{0^*\}$;
\item the transition probabilities are as follows: $p(0|0^*,1/n)=1-q_n,$ $p(n|0^*,1/n)=q_n,$ $p(0|0^*,0)=1$
and $p(0^{**}|s,a_0)=1$ for all $s\in S\setminus \{0^*\}$; the value of $q_n\in (0,1)$ for every $n\in \mathbb{N}$ will be
specified later;
\item the payoffs are as follows: $r(0^{**},a_0)=r(0^*,a)=0$ for every $a\in A(0^*)$ and $r(0,a_0)=1$, $r(s,a_0)=s$ for $s\in \mathbb{N}.$
\end{itemize}
Then, the function
$$w(s):=  s\mbox{ for }  s\in\mathbb{N}\quad\mbox{and}\quad w(s):=
1 \mbox{ for }s\in\{0,0^*,0^{**}\}
$$
is continuous on $S$ and satisfies (A1). Let $\mu=\delta_{0^*}$ and let $\pi_n$ denote the policy that chooses
the action $1/n$ in state $0^*.$
 We shall consider two cases.\\

(I) Let $q_n=1/n$ for $n\in\mathbb{N}.$ Note that $p$ is both weakly and setwise continuous.
We now check that (\ref{a21}) in (A2) does not hold.
Fix $l\in\mathbb{N}$ and choose $n>l$. Clearly, we have that
$$E_\mu^{\pi_n}(w(s_2)1_{\{w(s_2)\ge l\}})=1.$$ Then,
$$\sup_{n\in\mathbb{N}}E_\mu^{\pi_n}(w(s_2)1_{\{w(s_2)\ge l\}})=1$$
and consequently,
$$\lim_{l\to\infty}\sup_{\pi\in\Pi}E_\mu^{\pi}(w(s_2)1_{\{w(s_2)\ge l\}})=1,$$
so (\ref{a21}) does not hold.
On the other hand,
$$E_\mu^{\pi}w(s_k)=1\quad\mbox{ for any policy }\pi\in\Pi\mbox{ and }k>2.$$
Therefore,
$$\lim_{n\to\infty}\sup_{\pi\in\Pi} E_\mu^{\pi}\left(\sum_{k=n}^\infty \beta^k w(s_k)\right)=0.$$
Thus, (\ref{a12}) in (A2) holds.
Moreover, we observe that  (W4) and  (S4) are not met, since
$$\sum_{t\in S} w(t)p(t|0^*,1/n)=2-1/n \nrightarrow \sum_{t\in S} w(t)p(t|0^*,0)=1 \quad\mbox{as }n\to\infty.$$
Further, it is not difficult to see that (A3) holds with $\delta=2.$  In this example, the optimal policy does not exist, since
$$\sup_{n\in\mathbb{N}}E_\mu^{\pi_n}\left(\sum_{n=1}^\infty \beta^k r(s_n,a_n)\right)=\sup_{n\in\mathbb{N}}\beta(2-1/n)=2\beta$$
and for each policy $\pi_n$ the expected discounted payoff is strictly less than $2\beta.$  \\

(II) Let $q_n=1/2^n$ for $n\in\mathbb{N}.$ Obviously,  $p$ is again weakly and
setwise continuous. For these values of $q_n$
condition (\ref{a21}) in (A2) holds, because for $n>l$  we have that
$$\sup_{\pi\in\Pi}E_\mu^{\pi}(w(s_2)1_{\{w(s_2)\ge l\}})=l/2^l\quad\mbox{and}\quad
\lim_{l\to\infty}\sup_{\pi\in\Pi}E_\mu^{\pi}(w(s_2)1_{\{w(s_2)\ge l\}})=0.$$
Moreover, assumptions (W4) and  (S4) are also satisfied, because
$$\sum_{t\in S} w(t)p(t|0^*,1/n)=1-1/2^n+n/2^n\to \sum_{t\in S} w(t)p(t|0^*,0)=1 \quad\mbox{as }n\to\infty.$$
Observe that condition (\ref{a21}) is satisfied and  assumption  (A3) holds with $\delta=5/4.$
An optimal policy exists, because
$$\sup_{n\in\mathbb{N}}E_\mu^{\pi_n}\left(\sum_{n=1}^\infty \beta^k r(s_n,a_n)\right)=\sup_{n\in\mathbb{N}}\beta(1-1/2^n+n/2^n)=5\beta /4$$
and this supremum is realized by $\pi_2$ or $\pi_3.$ Obviously, we may replace the probabilities
$1/2^n,$ $n\in\mathbb{N},$
by arbitrary values of $q_n\in(0,1)$ such that $nq_n\to 0$ as $n\to\infty.$

%\begin{ack}                               % Place
% The authors acknowledge the financial support from the National Science Centre,
%Poland: Grant 2016/23/B/ST1/00425.
%\end{ack}

\bibliographystyle{plain}

\end{document}